# The Harmonic Variational Principle for the Einstein-Hilbert Functional

Sergey Stepanov, Irina Tsyganok

**Abstract.** Let $(M, g)$ be a compact $n$-dimensional Riemannian manifold, $n \geq 3$. We introduce a restricted variational principle for the Einstein–Hilbert functional by requiring the admissible metric variations to satisfy the harmonic gauge condition. We derive the corresponding Euler–Lagrange equation and show that a metric is critical with respect to all volume-preserving harmonic variations if and only if its Einstein tensor differs from a multiple of the metric by an element of the image of the adjoint Bianchi operator. We prove that every harmonic critical metric determines a compact Ricci soliton whose soliton constant is given by the normalized Einstein–Hilbert functional. By Perelman's theorem, every such metric is in fact the metric of a compact gradient Ricci soliton. Conversely, every compact gradient Ricci soliton satisfies the restricted Euler–Lagrange equation. Thus, a compact Riemannian metric is harmonic critical if and only if it is the metric of a compact gradient Ricci soliton. We further show that the gauge one-form differs from the negative differential of a soliton potential by a Killing one-form. In particular, if the Ricci tensor is negative definite, then the gauge one-form vanishes and the metric is Einstein. Moreover, every non-Einstein harmonic critical metric is necessarily shrinking.

## 1. Introduction

Let $M$ be a compact smooth manifold of dimension $n \geq 3$. The Einstein–Hilbert functional

$$\mathcal{S}(g) = \int_M s_g \, dv_g,$$

where $s_g$ and $dv_g$ denote the scalar curvature and the volume form of a Riemannian metric $g$, occupies a central position in Riemannian geometry and general relativity. A classical result of its variational theory states that Einstein metrics are precisely the critical points of $\mathcal{S} = \mathcal{S}(g)$ on the space of Riemannian metrics of fixed total volume. Equivalently, if the admissible infinitesimal variations range over all smooth symmetric 2-tensors satisfying the infinitesimal volume constraint, then the corresponding Euler–Lagrange equation is the Einstein equation. We refer to Berger and Ebin [2], Koiso [9], Schoen [15], and Besse [1] for the classical theory.

The present paper is motivated by the following natural question:

*To what extent can the class of admissible metric variations be restricted while preserving a meaningful geometric variational theory?*

Our answer to this question is provided by the *harmonic gauge*. More precisely, we will consider the infinite dimensional vector space of symmetric 2-tensors $h$ satisfying the natural first-order differential equation

$$B_g(h) := \delta h + \frac{1}{2} d\big(\mathrm{tr}_g h\big) = 0, \qquad (1.1)$$

where $B_g$ is the *Bianchi operator* acts on a symmetric 2-tensors (see [10]) and $\delta$ is the ordinary *divergence operator*. For example, the second *Bianchi identity* gives (see [1, p. 120])

$$B_g\big(\mathrm{Ric}_g\big) := \delta\, \mathrm{Ric}_g + \frac{1}{2}\, d\, s_g = 0$$

for the Ricci tensor $\mathrm{Ric}_g$ and the scalar curvature $s_g$ of an arbitrary Riemannian manifold $(M, g)$.

Condition (1.1) was introduced by Chen and Nagano in their theory of harmonic metrics and harmonic symmetric tensors [3, 4]. Following this terminology, symmetric 2-tensors satisfying (1.1) will be called *harmonic*, and the variations generated by them will be called *harmonic variations* throughout this paper.

**Remark 1.1.** Condition (1.1) coincides formally with the classical harmonic, or de Donder, gauge introduced in general relativity by de Donder [5] and later adapted to the Ricci flow by Hamilton [7] and DeTurck [20]. In geometric analysis, the same condition is commonly used as the linearized Bianchi gauge to eliminate diffeomorphism invariance; for a recent application to Ricci-flat conifolds, see Kröncke and Szabó [10]. Its role in the present paper is different: following Chen and Nagano [3, 4], we interpret $B_g(h) = 0$ as the infinitesimal harmonicity condition for the identity map and, consequently, as a restriction on the admissible directions in the variational theory of the Einstein–Hilbert functional.

To describe this equation, let

$$B_g^*\theta = \delta^*\theta + \frac{1}{2}(\delta\theta)g$$

be the formal $L^2$-adjoint of the Bianchi operator acts on 1-forms $\theta$. The Berger–Ebin decomposition (see [2]) associated with $B_g$ gives an $L^2$-orthogonal splitting

$$C^\infty(S^2T^*M) = \operatorname{Im} B_g^* \oplus \ker B_g\,.$$

Using this decomposition, we prove that a metric $g$ is critical with respect to all volume-preserving harmonic variations if and only if there exists a smooth 1-form $\theta$ such that

$$E_g = B_g^*\theta - \frac{n-2}{2n}\bar{s}_g\, g, \tag{1.2}$$

where

$$E_g = \operatorname{Ric}_g - \frac{1}{2}s_g g$$

is the standard *Einstein tensor* and

$$\bar{s}_g = \frac{1}{\operatorname{Vol}(M,g)}\int_M s_g\; dv_g$$

is the *average scalar curvature* of $(M, g)$. Metrics satisfying (1.2) will be called *harmonic critical metrics*.

The structural equation (1.2) has a particularly strong geometric interpretation. Taking its trace yields

$$\bar{s}_g = s_g + \delta\theta,$$

and therefore (1.2) reduces to

$$\operatorname{Ric}_g = \delta^*\theta + \frac{\bar{s}_g}{n}g. \tag{1.3}$$

If $X = \theta^\#$, then

$$\delta^*\theta = \frac{1}{2}\mathcal{L}_X g,$$

so that

$$\mathrm{Ric}_g + \frac{1}{2}\mathcal{L}_{-X} g = \lambda g, \qquad \lambda = \frac{\bar{s}_g}{n}. \tag{1.4}$$

Thus, every compact harmonic critical metric is the metric of a *Ricci soliton,* denoted by $(M. g, X)$ (see [6]). and its *soliton constant* $\lambda$ is determined by the normalized Einstein–Hilbert functional:

$$\lambda = \frac{1}{n\,\mathrm{Vol}(M, g)} \mathcal{S}(g). \tag{1.5}$$

When $X$ is the gradient vector field $\nabla f$ of a function $f$ on $(M, g)$ such a manifold is called *gradient Ricci soliton* $(M, g, \nabla f)$ with a potential $f$. In this case the previous equation becomes

$$Ric \; + \mathrm{Hess}_g f \; = \; \lambda g,$$

where $\mathrm{Hess}_g f$ is the standard Hessian of $f$. We say that a compact gradient soliton $(M. g, \nabla f)$ is *trivial* if its potential function $f$ is constant.

In turn, *Perelman's theorem* (see [6, Theorem 3.1]) asserting that every compact Ricci soliton is gradient, then implies that every compact harmonic critical metric is the metric of a gradient Ricci soliton $(M. g, \nabla f)$. Conversely, we show that every compact gradient Ricci soliton satisfies the restricted Euler–Lagrange equation (1.2). We therefore obtain the main characterization of the paper:

*A compact metric is harmonic critical if and only if it is the metric of a gradient Ricci soliton.*

Under this equivalence, the gauge 1-form $\theta$ is related to a soliton potential $f$ by

$$\theta = -df + \eta,$$

where $\eta^{\#}$ is a Killing vector field. Hence, in the absence of nonzero Killing fields, the harmonic gauge form is exactly the negative differential of the soliton potential.

The characterization also gives immediate rigidity consequences expressed directly in terms of the Einstein–Hilbert functional. If the average scalar curvature $\bar{s}_g = 0$, then the corresponding compact soliton is trivial, $g$ is Ricci-flat, and $\theta$ is a Killing 1-form. If the average scalar curvature $\bar{s}_g = 0$, then $g$ is Einstein with negative Einstein constant and $\theta = 0$. Consequently, every non-Einstein compact harmonic critical metric necessarily satisfies $\bar{s}_g > 0$.

The restricted variational principle developed here therefore provides a variational characterization of compact gradient Ricci solitons. Einstein metrics form its trivial branch, whereas its nontrivial branch consists of compact non-Einstein shrinking gradient Ricci solitons. In this way, the harmonic gauge connects the classical variational theory of the Einstein–Hilbert functional with the geometry of Ricci solitons.

The paper is organized as follows. Section 2 introduces the operators $B_g$and $B_g^*$, reviews harmonic symmetric 2-tensors, and establishes the required $L^2$-orthogonal decompositions. Section 3 develops the restricted variational principle and derives the Euler–Lagrange equation (1.2). Section 4 proves the equivalence between compact harmonic critical metrics and compact gradient Ricci solitons and derives the corresponding rigidity consequences.

**Remark**. An earlier version of part of the variational construction appeared in the preprint [22]. The present paper provides a substantially revised and self-contained treatment. In particular, the Ricci almost soliton interpretation of [22] is replaced by the exact equivalence with compact Ricci solitons; the converse implication and the Killing-field ambiguity are made explicit; and the analysis is supplemented by integral identities, rigidity results, and concrete non-Einstein examples.

## 2. Harmonic symmetric 2-tensors and their $L^2$-orthogonal decomposition

Let $(M, g)$ be a compact Riemannian manifold of dimension $n \geq 3$. Denote by $C^\infty(S^2M)$ the Fréchet space of smooth symmetric covariant 2-tensor fields on $M$, and by $C^\infty(T^*M)$ the space of smooth 1-forms.

For tensor $h, k \in C^\infty(S^2T^*M)$, their local components in a coordinate system $x^1, \dots, x^n$ will be denoted by $h_{ij} := h\left(\frac{\partial}{\partial x^i}, \frac{\partial}{\partial x^j}\right)$ and $k_{ij} := k\left(\frac{\partial}{\partial x^i}, \frac{\partial}{\partial x^j}\right)$ with respect to coordinate vector fields $\frac{\partial}{\partial x^j}$. The metric $g$ induces the natural $L^2$-inner product

$$\langle h, k\rangle_{L^2} = \int_M g^{ik}\, g^{jl} h_{ij} k_{kl}\, dv_g,$$

where $g_{ij} := g\left(\frac{\partial}{\partial x^i}, \frac{\partial}{\partial x^j}\right)$ and $\left(g^{ij}\right)$ are the contravariant components of the inverse metric tensor $g^{-1}$.

Let $\delta: C^\infty(S^2M) \to C^\infty(T^*M)$ be the divergence operator locally given by

$$\left(\delta h\right)_i = -\nabla^j h_{ji},$$

where $\nabla^i = g^{ij}\nabla_j$, and $\nabla_j$ denotes the covariant derivative with respect to the coordinate vector field $\frac{\partial}{\partial x^j}$. Then the classical Killing operator $\delta^*: C^\infty(T^*M) \to C^\infty(S^2M)$ be its $L^2$-formal adjoint operator defined by (see [2])

$$(\delta^*\theta)_{ij} = \frac{1}{2}(\nabla_i\theta_j + \nabla_j\theta_i).$$

where $\theta_j = \theta\left(\frac{\partial}{\partial x^j}\right)$ are the local covariant components of the 1-form $\theta$. Direct calculations show (see [2]) that $\delta^*$ has an injective principal symbol. In other words, $\delta^*$ is an *overdetermined elliptic operator* [1, Appendix].

Define the first-order differential operator $B_g^*: C^\infty(T^*M) \to C^\infty(S^2M)$ by

$$B_g^*(\theta) = \delta^*\theta + \frac{1}{2}(\delta\theta)g.$$

Let us consider the simplest properties of $B_g^*$. A direct calculation shows that the principal symbol of $B_g^*$ is injective (see Berger–Ebin [2, p. 385]). In addition, for $n \geq 3$, taking the trace of $B_g^*\theta = \delta^*\theta + \frac{1}{2}(\delta\theta)g$ gives

$$\mathrm{tr}_g(B_g^*\theta) = \frac{n-2}{2}\delta\theta.$$

Hence $B_g^*\theta = 0$ implies $\delta\theta = 0$, and therefore $\delta^*\theta = 0$. Thus, for $n \geq 3$, the kernel of $B_g^*$ coincides with the space of Killing 1-forms or, in other words, infinitesimal isometric transformations on $(M, g)$ (see [1, p. 40]).

The following identity will be useful later.

$$| B_g^*\theta |^2 = | \delta^*\theta |^2 - (\delta\theta)^2 + \frac{n}{4}(\delta\theta)^2 = | \delta^*\theta |^2 + \frac{n-4}{4}(\delta\theta)^2.$$

Therefore, if $n \geq 4$ (resp. $n \leq 4$), then

$$| B_g^*\theta |^2 \geq | \delta^*\theta |^2 \quad \left(\text{resp. } | B_g^*\theta |^2 \leq | \delta^*\theta |^2\right).$$

In particular, $| B_g^*\theta |^2 = | \delta^*\theta |^2$ for $n = 4$.

We now compute the formal $L^2$-adjoint of $B_g^*$. For this we consider arbitrary tensors

$$h \in C^\infty(S^2M), \qquad \theta \in C^\infty(T^*M).$$

Using integration by parts, we obtain

$$\langle B_g^*(\theta), h\rangle_{L^2} = \langle \delta^*\theta, h\rangle_{L^2} + \frac{1}{2}\int_M (\delta\theta)\,(\operatorname{tr}_g h)\,dv_g,$$

where $\delta\theta = -\nabla^i\theta_i$ is the divergence of the 1-form $\theta$.

Since

$$\langle \delta^*\theta, h\rangle_{L^2} = \langle \theta, \delta h\rangle_{L^2},$$

and

$$\langle \delta\theta, f\rangle_{L^2} = \langle \theta, df\rangle_{L^2}$$

for every smooth function $f \in C^\infty(M)$, we obtain

$$\frac{1}{2}\int_M (\delta\theta)(\operatorname{tr}_g h)\,dv_g = \frac{1}{2}\langle \theta, d(\operatorname{tr}_g h)\rangle_{L^2}.$$

Therefore,

$$\langle B_g^*(\theta), h\rangle_{L^2} = \langle \theta, \delta h + \frac{1}{2}d(\operatorname{tr}_g h)\rangle_{L^2}.$$

Hence, the formal $L^2$-adjoint of $B_g^*$ is the classical Bianchi operator

$$B_g(h) = \delta h + \frac{1}{2}d(\operatorname{tr}_g h).$$

Moreover, $\left(B_g^*\right)^* = B_g$, and $B_g^*$ is the formal $L^2$-adjoint of $B_g$.

**Remark 2.1.** The operator $B_g$coincides with the classical Bianchi operator used in modern geometric analysis. In the present paper, however, its kernel is interpreted as the natural space of admissible variations for the restricted Einstein–Hilbert variational principle.

Define the vector space of symmetric harmonic 2-tensors by

$$\mathcal{H}_g = \ker B_g = \left\{h \in C^\infty(S^2M): \delta h = -\frac{1}{2} d(\mathrm{tr}_g\, h)\right\}.$$

The operator $B_g$ is closely related to the classical Killing operator and to the conformal deformation operator appearing in the Berger–Ebin decomposition theory [2]. Since the injectivity of the principal symbol of $B_g^*$ is easily verified by a standard algebraic calculation in local coordinates, $B_g^*$ defines a first-order linear differential operator with injective principal symbol for $n \geq 2$ (see also [1, p. 456]; [2, p. 385]). Since the principal symbol of $B_g^*$ is injective, the second-order operator

$$\mathcal{A}_g = B_g B_g^*$$

is elliptic and formally self-adjoint (see [1, Appendix] and [2]). Hence, $\ker B_g B_g^* = \ker B_g^*$ and hence $\dim \ker \mathcal{A}_g < \infty$.

It is directly proved that the operator $B_g^*$ satisfies the identity

$$\left(\delta B_g^*\right)(\theta) = \tfrac{1}{2}\Delta_S \theta,$$

where

$$\Delta_S \theta := 2\delta\, \delta^* \theta - d(\delta\theta) \tag{2.1}$$

is the *Sampson Laplacian* defined on smooth 1-forms (see [18], [19]). In addition, $\Delta_S \theta = 0$ if and only if $\theta$ be the $g$-dual 1-form of an infinitesimal harmonic transformation $X$, i.e. $\theta = X^\flat$ (see [18]).

Using the above, one can prove that $\mathcal{A}_g := B_g B_g^*$ has the form

$$\mathcal{A}_g := B_g B_g^* = \frac{1}{2}\Delta_S + \frac{n-2}{4} d\delta.$$

By the classical Berger–Ebin–York decomposition theory based on Hodge–De Rham methods (see [1] and [2]), one obtains the $L^2$-orthogonal decomposition

$$C^\infty(S^2M) = \mathrm{Im}(B_g^*) \oplus_{L^2} \ker(B_g) = \mathrm{Im}(B_g^*) \oplus_{L^2} \mathcal{H}_g. \tag{2.2}$$

Consequently, every symmetric tensor $\varphi \in C^{\infty}(S^2M)$ admits a unique decomposition of the form

$$\varphi = B_g^*(\theta) + h,$$

where $h \in \mathcal{H}_g$ is a harmonic tensor.

Since the principal symbol of $B_g^*$ is injective, then the adjoint operator $B_g$ possesses an infinite-dimensional kernel (see Berger-Ebin [2]). Consequently, the space of harmonic tensors $\mathcal{H}_g = \ker B_g$ is infinite-dimensional (see also [3] and [4]).

According, to Chen and Nagano, a symmetric tensor field $h \in C^{\infty}(S^2M)$ satisfying $B_g(h) = 0$ is called *harmonic*. Correspondingly, variations generated by such tensors will be called *harmonic variations*. Therefore, the space $\mathcal{H}_g$ defines on $(M, g)$ a natural restricted class of symmetric harmonic variations of the metric $g$.

Thus, $\mathcal{H}_g = \ker B_g$ is the natural infinite-dimensional space of symmetric harmonic 2-tensors and determines the admissible variational directions in the restricted variational problem considered below.

An additional decomposition will be useful for relating the condition to *infinitesimal harmonic transformations*.

Let

$$\mathcal{TT}_g(M) = \{h \in C^{\infty}(S^2M) : \delta h = 0,\ \mathrm{tr}_g\, h = 0\}$$

denote the infinite dimensional space of transverse-traceless tensors (see [1, Theorem 12.30]). By the classical Berger-Ebin-York decomposition theory (see [2]), the algebraic sum

$$\mathrm{Im}\, \delta^* + C^{\infty}(M)\, g$$

is closed in $C^{\infty}(S^2M)$, and one has the $L^2$-orthogonal decomposition (see, e.g., [1, pp. 118-119]; [9])

$$C^{\infty}(S^2M) = (\mathrm{Im}\, \delta^* + C^{\infty}(M)\, g) \oplus_{L^2} \mathcal{TT}_g(M).$$

Accordingly, every tensor $h \in C^{\infty}(S^2M)$ can be represented in the form

$$h = \delta^*\theta + \lambda g + h^{TT},$$

where

$$\theta \in C^\infty(T^*M), \qquad \lambda \in C^\infty(M), \qquad h^{TT} \in \mathcal{TT}_g(M).$$

Although the pair $(\theta,\ \lambda)$ need not be uniquely determined, the tensor

$$\delta^*\theta + \lambda g$$

is uniquely determined by the above orthogonal decomposition.

**Proposition 2.1.** *Let $(M, g)$ be a compact Riemannian manifold of dimension $n \geq 3$. Suppose that*

$$h = \delta^*\theta + \lambda g + h^{TT},$$

*where $h \in C^\infty(S^2M)$ and $h^{TT} \in \mathcal{TT}_g(M)$. Then*

$$2B_g(h) = \Delta_S\theta + (n-2)\, d\lambda. \tag{2.2}$$

*Consequently, $h$ is harmonic if and only if*

$$\Delta_S\theta + (n-2)\, d\lambda = 0. \tag{2.3}$$

*In particular, any two of the following three assertions imply the remaining one:*

(i) *$h$ is a harmonic tensor;*

(ii) *$\theta^\#$ is an infinitesimal harmonic transformation;*

(iii) *$\lambda$ is constant.*

**Proof.** Since

$$\delta h^{TT} = 0, \qquad \mathrm{tr}_g\, h^{TT} = 0, \qquad \delta(\lambda g) = -\, d\lambda,$$

and

$$\mathrm{tr}_g(\delta^*\theta) = -\,\delta\theta,$$

we have

$$\delta h = \delta\delta^*\theta - d\lambda$$

and

$$d(\mathrm{tr}_g\, h) = -d(\delta\theta) + n\, d\lambda.$$

Therefore,

$$2B_g(h) \quad = 2\delta h + d\left(\operatorname{tr} h_g\right) = 2\delta\delta^*\theta - d(\delta\theta) + (n-2)d\lambda =$$

$$= \Delta_S\theta + (n-2)d\lambda.$$

This proves (2.2), while (2.3) follows from the definition

$$\mathcal{H}_g = \ker B_g.$$

The final assertion follows immediately from (2.2), since $n \geq 3$. □

**Corollary 2.2.** *Let $(M, g)$ be a compact Riemannian manifold with negative-definite Ricci tensor. Suppose that*

$$h = \delta^*\theta + \lambda g + h^{TT}, \quad h^{TT} \in \mathcal{TT}_g(M),$$

*is a symmetric harmonic 2-tensor and that $\lambda$ is constant. Then $\theta = 0$. Consequently,*

$$h = \lambda g + h^{TT}.$$

*If, in addition,* h *is infinitesimally volume-preserving, namely*

$$\int_M tr_g\, h\, dv_g = 0,$$

*then $\lambda = 0$, and hence $h = h^{TT}$.*

**Proof.** Since $h$ is a symmetric harmonic 2-tensor and $\lambda$ is constant, equation (2.3) gives

$$\Delta_S\theta = 0.$$

Thus, $\theta^\sharp$ is an infinitesimal harmonic transformation. The Bochner-Yano formula for the Sampson Laplacian yields (see [18])

$$0 = \int_M \langle \Delta_S\theta, \theta\rangle\, dv_g = \int_M \left(| \nabla\theta^\sharp |^2 - \operatorname{Ric}_g(\theta^\sharp, \theta^\sharp)\right)\, dv_g.$$

Since $\operatorname{Ric}_g$is negative-definite, it follows that $\theta^\sharp = 0$, and therefore $\theta = 0$. Hence

$$h = \lambda g + h^{TT}.$$

Finally,

$$\operatorname{tr}_g h = -\,\delta\theta + n\lambda.$$

Integrating over $M$ and using $\int_M \delta\theta \, d\, v_g = 0$, we obtain

$$\int_M \mathrm{tr}_g \, h \, dv_g = n \, \lambda \; \mathrm{Vol}(M, g).$$

Thus, the infinitesimal volume-preserving condition implies $\lambda = 0$, and consequently $h = h^{TT}$. □

## 3. The Restricted Harmonic Variational Principle

The purpose of this section is to develop a restricted variational theory for the Einstein-Hilbert functional based on harmonic variations. Unlike the classical Einstein variational principle, where admissible variations range over the entire Fréchet space $C^\infty(S^2M)$ of smooth symmetric tensor fields, we restrict the variational problem to the infinite-dimensional space $\mathcal{H}_g$ of symmetric harmonic 2-tensors introduced in Section 2.

This restriction may be viewed as an infinite-dimensional analogue of a constrained extremum problem in classical variational calculus. The ambient Fréchet manifold of Riemannian metrics remains unchanged, while only the admissible tangent directions are constrained by the harmonic gauge condition::

$$\delta \, h = -\frac{1}{2} d\big(\mathrm{tr}_g h\big),$$

for $h \in C^\infty(S^2M)$. As we shall see, this restriction nevertheless preserves a substantial part of the geometric content of the classical Einstein–Hilbert variational principle.

Let $S(g)$ denote the Einstein-Hilbert functional on the space of smooth Riemannian metrics on a compact Riemannian manifold $(M, g)$:

$$\mathcal{S}(g) = \int_M s_g \;dv_g$$

Let $g(t)$ be a smooth one-parameter variation of $g$, with $g(0) = g$, and set

$$h = \frac{d}{dt}\bigg|_{t=0} g(t).$$

The classical first variation formula (see Besse [1, §4.21]) is

$$\mathcal{S}'_g(h) = -\int_M \langle E_g, h\rangle dv_g, \tag{3.1}$$

where $E_g$ is the standard Einstein tensor:

$$E_g = Ric_g - \frac{1}{2}s_g g.$$

To capture the constraint of volume preservation, we require that the smooth variation $g(t)$ does not alter the total volume of the manifold $M$. To impose the infinitesimal volume constraint, we require

$$\left.\frac{d}{dt}\right|_{t=0} \mathrm{Vol}\big(M, g(t)\big) = 0.$$

Equivalently,

$$\int_M \mathrm{tr}_g\, h\, dv_g = 0. \tag{3.2}$$

Since the metric tensor $g$ itself belongs to the infinite-dimensional space $\mathcal{H}_g$ of harmonic tensors, condition (3.2) is equivalent to the global $L^2$-orthogonality relation $h \perp g$ for

$$g^{\perp} = \{h \in C^{\infty}(S^2T^*M) \colon \langle h, g\rangle_{L^2} = 0\}.$$

Consequently, the admissible volume-preserving harmonic variations are precisely characterized as symmetric tensor fields belonging to the intersection $h \in \mathcal{H}_g \cap g^{\perp}$. By the $L^2$-orthogonal decomposition established in Section 2, the Einstein tensor admits a decomposition

$$E_g = B_g^*(\theta) + E_{g.}^H. \tag{3.3}$$

where $\theta \in C^{\infty}(T^*M)$ and $E_g^{\mathcal{H}} \in \mathcal{H}_g$. We shall refer to $E_g^H$ as the harmonic component of the Einstein tensor.

The philosophy of the present construction differs essentially from other constrained variational approaches to the Einstein-Hilbert functional. In the theory of critical point metrics and related developments, one restricts the ambient space of metrics itself by imposing geometric conditions, such as constant scalar curvature. In contrast, the present approach preserves the full Fréchet manifold of Riemannian metrics and imposes a differential constraint only on the admissible deformation tensors. Thus, the restriction acts on tangent directions rather than on the underlying space of metrics.

The following theorem constitutes the main result of the paper.

**Theorem 3.1.** *Let $(M, g)$ be a compact Riemannian manifold of dimension $n \geq 3$. Then $g$ is a critical point of the Einstein-Hilbert functional with respect to all volume-preserving harmonic variations if and only if there exists a smooth 1-form $\theta$ such that*

$$E_g = B_g^*(\theta) - \frac{n-2}{2n} \bar{s}\, g.$$

*Moreover, $\theta^\#$ is an infinitesimal harmonic transformation.*

**Proof.** Assume that $g$ is critical with respect to all admissible harmonic variations, meaning that the first variation (3.1) vanishes for all $\varphi \in \mathcal{H}_g \cap g^\perp$. Using

decomposition (3.3), write:

$$E_g = B_g^*(\theta) + E_{g.}^H.$$

Substituting this decomposition into (3.1), we obtain:

$$\int_M \langle E_g, h\rangle \, dv_g = 0$$

for every $h \in \mathcal{H}_g \cap g^\perp$. Since $\operatorname{Im} B_g^*$ is $L^2$-orthogonal to $\mathcal{H}_g$ by the properties of the Berger–Ebin decomposition (see [2]), it follows that:

$$\int_M \langle E_{g.}^H, h\rangle dv_g = 0$$

for every $h \in \mathcal{H}_g \cap g^\perp$. Hence, we have the geometric orthogonality:

$$E_g^H \perp \left( \mathcal{H}_g \cap g^\perp\right).$$

Since $g \in \mathcal{H}_g$, the Hilbert-space decomposition

$$\mathcal{H}_g = \left(\mathcal{H}_g \cap g^\perp\right) \oplus \operatorname{span}\{g\}$$

holds. Therefore,

$$E_g^{\mathcal{H}} = cg$$

for some constant $c \in \mathbb{R}$.

Substituting this relation back into (3.3), we obtain the structural balance equation:

$$E_g = B_g^*(\theta) + c\, g \tag{3.4}$$

Conversely, assume that equation (3.4) holds. Let $h \in \mathcal{H}_g \cap g^{\perp}$ be an arbitrary admissible harmonic variation. By the orthogonality of the Berger–Ebin decomposition, we have:

$$\int_M \langle B_g^* (\theta) , h \rangle \, dv_g = 0$$

while the volume-preserving condition yields:

$$\int_M \langle cg , h \rangle \, dv_g = 0$$

because $h \perp g$. Consequently, summing these equations gives:

$$\int_M \langle E_g , h \rangle \, dv_g = 0$$

and the classical variation formula (3.1) implies $\frac{d}{dt}|_{t=0}\, \mathcal{S}(g(t)) = 0.$

Hence, $g$ is a genuine critical point of the Einstein-Hilbert functional with respect to all volume-preserving harmonic variations.

Taking the trace of equation (3.4), we obtain $s_g = -\,\delta\theta - \frac{2n}{n-2} c$. Integrating over the compact manifold $M$ and using the divergence theorem, we have

$$\int_M \delta\theta \, dv_g = 0.$$

Therefore, the total scalar curvature satisfies

$$\mathcal{S}(g) = \int_M s_g \; dv_g = -\frac{2n}{n-2} c \,\mathrm{Vol}(M, g),$$

Consequently,

$$c = -\frac{n-2}{2n} \cdot \frac{\mathcal{S}(g)}{\mathrm{Vol}(M, g)},$$

where

$$\bar{s}_g = \frac{1}{\mathrm{Vol}(M,g)} \int_M s_g \, dv_g$$

is the *average scalar curvature* over $(M, g)$ (see, for example, [1]; [7]; [20]). Then

$$c = -\frac{n-2}{2n} \cdot \bar{s}_g.$$

Thus, the constant c appearing in the structural equation is completely determined by the average scalar curvature of the metric and, in the unit-volume case, by the total scalar curvature.

Applying the operator $\delta$ to both parts of equality (3.4), we obtain $\Delta_S \theta = 0$ since $\delta E_g = 0$ and $\delta(cg) = 0$. Then, from (3.4), we deduce $\left(\delta B_g^*\right)(\theta) = 0$. At the same time, using the identity $\left(\delta B_g^*\right)(\theta) = \frac{1}{2}\Delta_S \theta$, we obtain $\Delta_S \theta = 0$. Therefore, if $\theta^\# = X$ is the $g$-dual vector field of $\theta$, then $X$ is an infinitesimal harmonic transformation. Theorem is proved. □

**Definition 3.1.** A Riemannian metric satisfying the restricted Euler–Lagrange equation (3.4) will be called harmonic critical.

**Remark 3.1.** Theorem 3.1 shows that restricting the admissible tangent directions by the harmonic gauge still leads to a genuine Euler–Lagrange equation governed by the Einstein tensor. Unlike variational theories in which the ambient space of metrics is restricted, the present construction retains the full Fréchet manifold of Riemannian metrics and imposes the constraint only on its tangent directions.

**Remark 3.2.** In the gradient case $\theta = df$, equation (3.4) takes the form

$$E_g = \left(\mathrm{Hess}\, f - \frac{1}{2}(\Delta f) g\right) + c\, g.$$

Consequently, the restricted Euler–Lagrange equation exhibits a formal similarity to several classical geometric equations involving Hessian potentials, including those arising in the theory of critical point metrics and in the Fischer–Marsden theory of the scalar-curvature operator; see Besse [1, Proposition 4.47 and Remark 4.48]. Although the geometric origins of these equations are different, this resemblance suggests that harmonic critical metrics belong naturally to the broader family of variational structures governed by Hessian-type deformations of the Ricci tensor.

The terminology is justified by Theorem 3.1, which shows that equation (3.4) is precisely the Euler–Lagrange condition associated with the Einstein-Hilbert functional under volume-preserving harmonic variations. Thus, harmonic critical metrics play, within the restricted variational framework developed here, the same role that Einstein metrics play in the classical Einstein variational principle.

As shown in the next section, equation (3.4) implies that every compact harmonic critical metric is the metric of a Ricci soliton and hence, by Perelman's theorem, of a gradient Ricci soliton. Consequently, harmonic critical metrics provide a natural bridge between constrained variational geometry, orthogonal gauge decompositions, and modern Ricci flow theory.

## 4. Harmonic Critical Metrics and Gradient Ricci Solitons

The purpose of this section is to determine the geometric meaning of the restricted Euler–Lagrange equation obtained in Theorem 3.1. We first show that every harmonic critical metric determines a Ricci soliton whose soliton constant is explicitly expressed in terms of the normalized Einstein–Hilbert functional. We then apply Perelman's theorem asserting that every compact Ricci soliton is gradient. Finally, we prove the converse and obtain a variational characterization of compact gradient Ricci solitons.

Recall that a harmonic critical metric satisfies

$$E_g = B_g^* \theta - \frac{n-2}{2n} \, \bar{s}_g \, g, \tag{4.1}$$

where

$$E_g = \mathrm{Ric}_g - \frac{1}{2} s_g g, \qquad B_g^* \theta = \delta^* \theta + \frac{1}{2} (\delta \theta) g,$$

and

$$\bar{s}_g = \frac{1}{\mathrm{Vol}(M,g)} \int_M s_g \, dv_g$$

is the average scalar curvature of $(M, g)$.

Substituting these expressions into (4.1), we obtain

$$\mathrm{Ric}_g - \frac{1}{2} s_g g = \delta^* \theta + \frac{1}{2} (\delta\theta) g - \frac{n-2}{2n} \bar{s}_g g. \tag{4.2}$$

Hence,

$$\mathrm{Ric}_g = \delta^* \theta + \left( \frac{1}{2} s_g + \frac{1}{2} \delta\theta - \frac{n-2}{2n} \bar{s}_g \right) g. \tag{4.3}$$

Taking the trace of (4.1) and using $\mathrm{tr}_g E_g = -\frac{n-2}{2} s_g$ and $\mathrm{tr}_g\left(B_g^* \theta\right) = \frac{n-2}{2} \delta\theta$, we obtain

$$-\frac{n-2}{2} s_g = \frac{n-2}{2} \delta\theta - \frac{n-2}{2} \bar{s}_g.$$

Since $n \geq 3$, it follows that

$$s_g + \delta\theta = \bar{s}_g. \tag{4.4}$$

Equivalently,

$$\delta\theta = \bar{s}_g - s_g. \tag{4.5}$$

Let $X = \theta^\#$. Since $\delta\theta = -\mathrm{div}_g X$, equation (4.5) becomes

$$\mathrm{div}_g X = s_g - \bar{s}_g. \tag{4.6}$$

Thus, the divergence of $X$ measures precisely the deviation of the scalar curvature from its average value.

Substituting (4.4) into (4.3), we obtain

$$\mathrm{Ric}_g = \delta^* \theta + \frac{\bar{s}_g}{n} g. \tag{4.7}$$

Since $\delta^* \theta = \frac{1}{2} \mathcal{L}_X g$, equation (4.7) can be written as

$$\mathrm{Ric}_g = \frac{1}{2}\mathcal{L}_X g + \lambda g,\ \lambda = \frac{\bar{s}_g}{n}. \tag{4.8}$$

Equivalently,

$$\mathrm{Ric}_g + \frac{1}{2}\mathcal{L}_{(-X)} g = \lambda g. \tag{4.9}$$

Thus, $g$ is a Ricci soliton with soliton vector field $-X$. Notice that $X = \theta^{\#}$ is the vector field naturally associated with the gauge form, whereas the vector field occurring in the standard Ricci soliton equation is $-X$.

We have proved the following result.

**Theorem 4.1.** Let $(M, g)$ be a compact Riemannian manifold of dimension $n \geq 3$. If $g$ is a harmonic critical metric, then $g$ is the metric of a Ricci soliton. More precisely, if $\theta$ is the gauge one-form occurring in the restricted Euler–Lagrange equation and $X = \theta^{\#}$, then

$$\mathrm{Ric}_g + \frac{1}{2}\mathcal{L}_{(-X)} g = \lambda g,$$

where

$$\lambda = \frac{\bar{s}_g}{n} = \frac{1}{n\mathrm{Vol}(M, g)} \int_M s_g \, dv_g. \tag{4.10}$$

Moreover,

$$\mathrm{div}_g X = s_g - \bar{s}_g. \tag{4.11}$$

**Remark 4.1.** A Ricci soliton is called shrinking, steady, or expanding according as

$$\lambda > 0, \qquad \lambda = 0, \quad \lambda < 0,$$

respectively; see [6]. Formula (4.10) shows that the type of the Ricci soliton associated with a harmonic critical metric is completely determined by the sign of the Einstein–Hilbert functional. More precisely, it is shrinking, steady, or expanding according as

$$\int_M s_g \, dv_g > 0, \quad \int_M s_g \, dv_g = 0, \quad \int_M s_g \, dv_g < 0,$$

respectively. On the space of unit-volume metrics, $\lambda = \frac{1}{n}\mathcal{S}(g)$.

**Remark 4.2.** If the restricted Euler–Lagrange equation is written in the form

$$E_g = B_g^*\theta + cg,$$

then taking its trace and integrating over $M$ gives

$$c = -\frac{n-2}{2n}\bar{s}_g. \tag{4.12}$$

Consequently,

$$\lambda = -\frac{2c}{n-2} = \frac{\bar{s}_g}{n}. \tag{4.13}$$

Thus, both constants occurring in the restricted Euler–Lagrange equation have natural variational interpretations: $c$is determined by the average scalar curvature, whereas the Ricci soliton constant $\lambda$ is given by the normalized Einstein–Hilbert functional.

We now use the following fundamental theorem of Perelman: every compact Ricci soliton is a gradient Ricci soliton; see [6, Theorem 3.1]. Combining this theorem with Theorem 4.1, we obtain the principal geometric consequence of the restricted harmonic variational principle.

**Theorem 4.2.** Let $(M, g)$ be a compact Riemannian manifold of dimension $n \geq 3$, and suppose that $g$is harmonic critical. Then there exists a smooth function $f \in C^\infty(M)$ such that

$$\mathrm{Ric}_g + \mathrm{Hess}_g f = \lambda g, \tag{4.14}$$

where

$$\lambda = \frac{\bar{s}_g}{n} = \frac{1}{n\mathrm{Vol}(M, g)}\int_M s_g \, dv_g. \tag{4.15}$$

Thus, every compact harmonic critical metric is the metric of a gradient Ricci soliton.

**Proof.** By Theorem 4.1, $g$ is the metric of a compact Ricci soliton. Perelman's theorem therefore implies the existence of a smooth function $f$ satisfying (4.14). Formula (4.15) follows directly from Theorem 4.1. □

The gradient representation supplied by Perelman need not coincide exactly with the original gauge representation. Indeed, equations (4.7) and (4.14) give

$$\delta^*\theta + \mathrm{Hess}_g f = 0.$$

Since $\mathrm{Hess}_g f = \delta^*(df)$, we obtain

$$\delta^*(\theta + df) = 0. \tag{4.16}$$

Consequently, the vector field $K = (\theta + df)^\#$ is Killing. Therefore,

$$\theta = -df + \eta, \tag{4.17}$$

where $\eta = K^\flat$ is a Killing one-form.

The identity (4.17) should be distinguished from the formal substitution $\theta = df$ considered in Remark 3.2. There, $\theta = df$ is merely a particular exact gauge ansatz used to illustrate the form of the restricted Euler–Lagrange equation. In contrast, equation (4.17) is an intrinsic consequence of Perelman's theorem: the function $f$ is a Ricci soliton potential, and the difference between the original gauge form and $-df$ is necessarily a Killing one-form.

We therefore obtain the following refinement.

**Corollary 4.3.** *Let $g$ be a harmonic critical metric on a compact manifold, let $\theta$ be the gauge one-form occurring in the restricted Euler–Lagrange equation, and let $f$ be a Perelman potential satisfying* (4.14). *Then*

$$\theta = -df + \eta,$$

*where $\eta^\#$ is a Killing vector field.*

In particular, if $(M, g)$ admits no nonzero Killing vector fields, then

$$\theta = -\, df. \tag{4.18}$$

Taking the trace of the gradient Ricci soliton equation (4.14) and using the convention $\mathrm{tr}_g\big(\mathrm{Hess}_g f\big) = -\,\Delta f$, we obtain $s_g - \Delta f = n\lambda$. By (4.15),

$$n\,\lambda = \bar{s}_g,$$

and hence

$$\Delta f = s_g - \bar{s}_g. \tag{4.19}$$

Thus, the potential function measures the deviation of the scalar curvature from its mean value.

Standard identities for gradient Ricci solitons give

$$ds_g = 2\mathrm{Ric}_g(\nabla, f\ \cdot) \tag{4.20}$$

and

$$s_g + |\,\nabla f\,|^2 - 2\lambda f = C \tag{4.21}$$

for some constant $C$; see [6, Proposition 2.1].

Combining these identities with Theorem 4.2, we obtain the following summary.

**Corollary 4.4.** *Let $(M, g)$ be a compact Riemannian manifold with a harmonic critical metric. Then there exists a smooth function $f \in C^\infty(M)$ such that*

$$\mathrm{Ric}_g + \mathrm{Hess}_g f = \frac{\bar{s}_g}{n} g,$$

*and* $\Delta f = s_g - \bar{s}_g$, $\quad ds_g = 2\mathrm{Ric}_g(\nabla, f\ \cdot)\ (\nabla f,\ \cdot)$, $\quad s_g + |\,\nabla f\,|^2 - 2\frac{\bar{s}_g}{n} f = C$ *for some constant* $C$.

The following consequences express the compact rigidity of steady and expanding Ricci solitons directly in terms of the Einstein–Hilbert functional.

**Corollary 4.5.** *Let $(M, g)$ be a compact Riemannian manifold with a harmonic critical metric. If*

$$\bar{s}_g = 0, \tag{4.22}$$

*or, equivalently,*

$$\int_M s_g \, dv_g = 0,$$

*then $g$ is Ricci-flat. Moreover, the gauge one-form $\theta$ is a Killing one-form.*

**Proof.** By Theorem 4.2, $g$ is the metric of a compact gradient Ricci soliton with soliton constant

$$\lambda = \frac{\bar{s}_g}{n}.$$

Condition (4.22) implies $\lambda = 0$. Every compact steady Ricci soliton is trivial; see [6, Proposition 3.4 and Corollary 3.5]. Hence $f$ is constant and $\mathrm{Ric}_g = 0$. By Corollary 4.3,

$$\theta = -df + \eta,$$

where $\eta^\#$ is Killing. Since $f$ is constant, $\theta = \eta$, and therefore $\theta$ is a Killing one-form. □

**Corollary 4.6.** *Let $(M, g)$ be a compact Riemannian manifold with a harmonic critical metric. If*

$$\bar{s}_g < 0, \tag{4.23}$$

*or, equivalently*,

$$\int_M s_g \, dv_g < 0,$$

*then $g$ is Einstein with negative Einstein constant and $\theta = 0$. More precisely*,

$$\mathrm{Ric}_g = \frac{\bar{s}_g}{n} g. \tag{4.24}$$

**Proof.** By Theorem 4.2, $g$ is the metric of a compact gradient Ricci soliton with

$$\lambda = \frac{\bar{s}_g}{n} < 0.$$

Every compact expanding Ricci soliton is trivial; see [6, Proposition 3.4 and Corollary 3.5]. Consequently, $f$ is constant and

$$\mathrm{Ric}_g = \lambda g = \frac{\bar{s}_g}{n} g.$$

In particular, the Ricci tensor is negative definite. By Corollary 4.3,

$$\theta = -df + \eta = \eta,$$

where $\eta^{\#}$ is Killing. Since a compact Riemannian manifold with negative-definite Ricci tensor admits no nonzero Killing vector fields, it follows that $\theta = 0$. □

**Remark 4.3.** Corollaries 4.5 and 4.6 show that every compact harmonic critical metric with nonpositive Einstein–Hilbert functional is Einstein. The conclusions concerning the gauge one-form are different, however. If

$$\int_M s_g \, dv_g = 0,$$

then $g$ is Ricci-flat and $\theta$ is Killing, but $\theta$ need not vanish. If

$$\int_M s_g \, dv_g < 0,$$

then $g$ is Einstein with negative Einstein constant and necessarily $\theta = 0$. The soliton potential also gives the following elementary rigidity criterion.

**Corollary 4.7.** *Let $(M, g)$ be a compact Riemannian manifold with a harmonic critical metric, and let $f$ be a potential function satisfying*

$$\mathrm{Ric}_g + \mathrm{Hess}_g f = \frac{\bar{s}_g}{n} g.$$

*If* $\langle \nabla s, \nabla f \rangle \leq 0$ *everywhere on M, then f and* $s_g$ *are constant, and* $g$ *is Einstein.*

**Proof.** By (4.19),

$$\Delta f = s_g - \bar{s}_g.$$

Using Green's formula, we obtain

$$\int_M \langle \nabla, s_g \nabla f \rangle \, dv_g = -\int_M s_g \, \Delta f \, dv_g.$$

Since

$$\int_M \bar{s}_g \left(s_g - \bar{s}_g\right) dv_g = 0,$$

it follows that

$$-\int_M s_g \, \Delta f \, dv_g = -\int_M (s_g - \bar{s}_g)^2 \, dv_g.$$

Therefore,

$$\int_M \langle \nabla, s_g \nabla f \rangle \, dv_g = -\int_M (s_g - \bar{s}_g)^2 \, dv_g.$$

Under the assumed inequality, both sides can vanish only if $s_g = \bar{s}_g$. Hence $\Delta f = 0$. Since $M$ is compact, $f$ is constant. The gradient Ricci soliton equation then reduces to $\mathrm{Ric}_g = \frac{\bar{s}_g}{n} g$, and $g$ is Einstein. □

Finally, the converse of Theorem 4.2 also holds.

**Theorem 4.8.** *Let* $(M, g, f)$ *be a compact gradient Ricci soliton of dimension* $n \geq 3$*, satisfying*

$$\mathrm{Ric}_g + \mathrm{Hess}_g f = \lambda g. \tag{4.25}$$

*Then* $\lambda = \frac{\bar{s}_g}{n}$*, and* $g$ *is a harmonic critical metric.*

**Proof.** Taking the trace of (4.25), we obtain

$$s_g - \Delta f = n\lambda. \tag{4.26}$$

Integrating over the compact manifold $M$ and using $\int_M \Delta f \, dv_g = 0$, we obtain

$$\lambda = \frac{1}{n\mathrm{Vol}(M,g)} \int_M s_g \, dv_g = \frac{\bar{s}_g}{n}. \tag{4.27}$$

Set $\theta = -df$. Since $\delta(df) = \Delta f$, we have $\delta\theta = -\Delta f$ and

$$B_g^* \theta = -\mathrm{Hess}_g f - \frac{1}{2}(\Delta f) g.$$

Using (4.25) and (4.26), we compute

$$\begin{aligned} E_g &= \mathrm{Ric}_g - \frac{1}{2} s_g g \\ &= -\mathrm{Hess}_g f + \left(\lambda - \frac{1}{2} s_g\right) g \\ &= -\mathrm{Hess}_g f - \frac{1}{2}(\Delta f) g - \frac{n-2}{2} \lambda g \\ &= B_g^* \theta - \frac{n-2}{2n} \bar{s}_g \, g. \end{aligned}$$

Thus, $g$ satisfies the restricted Euler–Lagrange equation of Theorem 3.1 and is therefore harmonic critical. □

Combining Theorems 4.2 and 4.8, we arrive at the main characterization.

**Corollary 4.9.** *A compact Riemannian metric $g$ of dimension $n \geq 3$ is harmonic critical if and only if it is the metric of a gradient Ricci soliton.*

Under this equivalence, the soliton constant is determined by the Einstein–Hilbert functional through

$$\lambda = \frac{1}{n\mathrm{Vol}(M,g)} \mathcal{S}(g).$$

Thus, the restricted harmonic variational principle provides a variational characterization of compact gradient Ricci solitons. Einstein metrics constitute the trivial branch of this variational theory. Moreover, Corollaries 4.5 and 4.6 show that every non-Einstein compact harmonic critical metric necessarily satisfies

$$\int_M s_g \, dv_g > 0$$

and therefore belongs to the shrinking gradient Ricci soliton branch.

## Conclusion and Further Directions

The results obtained above show that the harmonic gauge defines a genuinely geometric restricted variational principle for the Einstein–Hilbert functional. Its critical points are not merely related to Ricci solitons: on compact manifolds they coincide precisely with gradient Ricci solitons. Thus, the harmonic constraint produces a variational theory whose trivial branch consists of Einstein metrics and whose nontrivial branch consists necessarily of shrinking gradient Ricci solitons. The study of the second variation directly on the full space of volume-preserving harmonic variations, without reduction to the transverse-traceless gauge, remains a natural problem for future investigation.

## Data Availability Statement

No datasets were generated or analyzed during the current study.

## Funding Declaration

The author declares that no financial funds, grants or other support were used in the preparation of this manuscript.

## Conflict of interest statement

The author declares no conflict of interest.

Sergey Stepanov

Department of Mathematics All Russian Institute for Scientific and Technical Information of the Russia Academy of Sciences
20, Usievicha street, Moscow, Russia 125190

Department of Mathematics and Data Analysis, Finance University
49-55, Leningradsky Prospect, Moscow, Russia 125468
e-mail: s.e.stepanov@mail.ru

Irina Tsyganok

Department of Mathematics and Data Analysis, Finance University
49-55, Leningradsky Prospect Moscow Russia 125468
e-mail: i.i.tsyganok@mail.ru